\documentclass[11pt]{article}
\usepackage{amssymb,latexsym,amsmath,amsthm,mathrsfs,hyperref}
\newtheorem{theorem}{Theorem}[section]
\newtheorem{proposition}[theorem]{Proposition}
\newtheorem{lemma}[theorem]{Lemma}

\theoremstyle{remark} 
\newtheorem{remark}[theorem]{Remark} 

\theoremstyle{definition}

\def\<{\langle} 
\def\>{\rangle} 
\def\L{\mathcal{L}}
\def\K{\mathcal{K}}
\def\R{\mathbb{R}}
\def\C{\mathbb{C}}
\def\H{\mathbb{H}}
\def\S{\mathbb{S}}
\def\P{\mathcal{P}}
\def\Re{{\rm Re}}
\def\Im{{\rm Im}}
\def\zb{\bar z}
\def\d{{\rm d}}
\newcommand{\abs}[1]{|#1|}
\title{
Willmore Surfaces of Constant M\"obius Curvature
}
\author{Xiang Ma, Changping Wang
\footnote{This work is partially supported 
by RFDP (No. 20040001034)}}
\begin{document}
\maketitle

\begin{abstract}
We study Willmore surfaces of constant M\"obius curvature $\K$ 
in $\S^4$. It is proved that such a surface in $\S^3$ must be 
part of a minimal surface in $\R^3$ or the Clifford torus.
Another result in this paper is that
an isotropic surface (hence also Willmore) in $\S^4$ 
of constant $\K$ could only be part of a complex curve in 
$\C^2\cong \R^4$ or the Veronese 2-sphere in $\S^4$. 
It is conjectured that they are the only examples possible.
The main ingredients of the proofs are over-determined systems
and isoparametric functions.
\end{abstract}

\section{Introduction}

In this paper we consider Willmore surfaces in $\S^4$ of 
constant M\"obius curvature. 

Our interest in such surfaces 
originates from the analogous problem about minimal surfaces 
in $\S^n$. It was well known that minimal surfaces are quite 
rigid objects. In particular, minimal surfaces in $\S^n$ 
whose induced metric are of constant Gaussian curvature $K$
can be determined even locally (see \cite{bryant2} and 
references therein). Such a surface is always homogeneous, 
and locally it must be congruent to a generalized Veronese 
2-sphere (whose curvature is $2/m(m+1)$ for some integer $m$) 
when $K>0$, or to a generalized Clifford surface when $K=0$
(c.f. \cite{calabi,kenmotsu1} for construction of these examples); 
the case $K<0$ is impossible.
This remarkable result has been generalized to minimal surfaces 
in real space forms $\R^n, \H^n$ \cite{bryant2}
or in complex space form $\C\mathbb{P}^n$ 
(see \cite{kenmotsu3} and references therein).

In M\"obius geometry, associated with any surface without 
umbilic points there is a conformally invariant metric
(induced from the conformal Gauss map \cite{bryant1}) 
called the M\"obius metric, whose curvature 
is named the M\"obius curvature, and denoted by $\K$. 
The area of this surface under the M\"obius metric is exactly 
the well-known Willmore functional. Thus Willmore surfaces, 
the critical surfaces with respect to Willmore functional, 
are also called M\"obius minimal surfaces. 
More interestingly, minimal surfaces in space forms 
$\S^n,\R^n,\H^n$ are themselves also 
Willmore surfaces \cite{blaschke, wang}.
This raises the natural problem that whether we can 
generalize the remarkable classification result above to 
Willmore surfaces of constant M\"obius curvature. 

In 3-space the only known examples of such surfaces 
are minimal surfaces in $\R^3$ (with $\K\equiv 1$) and 
the Clifford torus (with $\K\equiv 0$). In 4-space we might
add complex curves in $\C^2\cong \R^4$ (with $\K\equiv 2$) 
and the Veronese sphere (with $\K\equiv \frac 12$). 
Naturally we have the following

\noindent {\bf Conjecture.~~}{\it
Any Willmore surface of constant M\"obius curvature is locally 
M\"obius equivalent to one of the four examples mentioned above
(minimal surfaces in $\R^3$, Clifford torus, complex curve
in $\C^2$, and Veronese sphere).}

We can confirm this conjecture under various additional 
assumptions. First comes the codim-1 case:

\noindent {\bf Theorem A.~~}{\it
Let $y:M\to \S^3$ be a Willmore surface without umbilic points 
and $\K=constant$. Then locally $y$ is M\"obius equivalent 
either to a minimal surface in $\R^3$ with $\K=1$,
or to the Clifford torus with $\K=0$.}

The second result concerns isotropic surfaces
in $\S^4$ (also known as super-conformal surfaces), which is an
important class of Willmore surfaces.

\noindent{\bf Theorem B.~~}{\it
Let $y:M\to \S^4$ be an isotropic surface without umbilic 
points and $\K=constant$. Then locally $y$ is M\"obius 
equivalent either to a complex curve in $\C^2$ 
with $\K=2$, or to the Veronese sphere with $\K=\frac 12$.}

These results are remarkable in several aspects. 
First, Theorem A gives a local characterization of 
the Clifford torus in the category of M\"obius geometry, 
which might be helpful in the exploration 
of Willmore Conjecture and other problems. 

Next, this is the starting point to the 
classification of all such surfaces in $\S^n$. 
Surely this problem is more difficult than the original 
one about minimal surfaces, whereas the importance of 
the latter problem and its solutions have already been recognized.

Finally, although M\"obius metric as well as M\"obius curvature 
are basic invariants in the surface theory of M\"obius geometry, 
so far they receive little attention, and there exists 
few results related to them (one exception is \cite{jeromin}). 
Our theorem might be the first deep result in this direction. 
In particular, it partly answers the following question: {\it 
To what extent is the M\"obius metric restricted 
when we assume the surface to be Willmore?}
Solution to this problem will help us to have 
a better understanding of Willmore surfaces. In general,
one may expect to study the \emph{intrinsic geometry} for
various interesting surface classes in M\"obius geometry.

Note that minimal surfaces in $\R^n$ or $\S^n$ are 
special examples of Willmore surfaces. Instead of studying
such surfaces with constant Gaussian curvature, we might 
consider them under the condition of constant $K_{III}$,
where $K_{III}$ is the curvature of the metric 
(the third fundamental form) induced from the generalized 
Gauss map \cite{fujiki,hoffman,sasaki1,sasaki2,sasaki3,vlachos1}.
This subject is included in our research program,
because now the generalized Gauss map coincides with the 
conformal Gauss map in M\"obius geometry, and $K_{III}=\K$.
Therefore, those old results gain new meaning from 
our viewpoint, and provide useful informations to our study.
(Note they are also related to the analogous subject:
holomorphic curves in $\C\mathbb{P}^n$ with constant curvature,
see \cite{lawson}.) In particular, as partial confirmation 
to the conjecture before, there are the following known results.

\noindent{\bf Theorem.(\cite{hoffman})~}
{\it Any minimal surface in $\R^4$ with constant $K_{III}(=\K)$
is locally congruent to a minimal surface in $\R^3$ or a 
complex curve in $\C^2$ up to rigid motion.} 

\noindent{\bf Theorem.(\cite{sasaki1,sasaki2})~}
{\it Any compact minimal surface in $\S^4$ with constant 
$K_{III}(=\K)$ and with isolated umbilic points
is congruent to the Veronese sphere in $\S^4$ 
when it is a minimal 2-sphere, or to the Clifford torus 
when it is a minimal torus in $\S^3$, up to rigid motion.} 

Surprisingly, all known examples of Willmore surfaces 
with constant $\K$ are minimal surfaces in space forms.
We guess that there exists no other examples. Yet it seems
very hard to prove it.

Theorem~A was partly obtained in \cite{gorokh2}, 
where the subject was studied under the name of 
{\it minimal surfaces of constant Gaussian curvature
in a pseudo-Riemannian sphere}. For related work see
\cite{gorokh1}. Theorem~A has also been generalized
to space-like Willmore surfaces in 3-dim conformal
Lorentzian geometry \cite{nie}.

It is observed that the theorems above are local results 
except the last one. The authors believe that global
assumptions are not necessary in such problems.

\section{Surface theory in M\"obius geometry}

In this part we will briefly review the surfaces theory of 
M\"obius geometry. The basic reference is \cite{bpp}.

As usual, the unit sphere $\S^n\subset \R^{n+1}$ is 
identified with the projectivized light cone via
\[\S^n\cong \mathbb{P}(\L):~~x\leftrightarrow [1,x],\]
where $\L\subset \R^{n+1,1}$ is the null cone in the
$(n+2)$-dim Minkowski space for the quadratic form
$\<Y,Y\>=-Y_0^2+\sum_{i=1}^{n+1}Y_i^2$. In this way,
a point in $\S^n$ corresponds to an isotropic line, and
a $k$-dim sphere is viewed as $\mathbb{P}(\L\cap U^{\bot})$, 
which is determined by a space-like $(n-k)$-dim subspace $U$
(or its orthogonal complement $V=U^{\bot}$). 
The projective action of the Lorentz group on 
$\mathbb{P}(\L)$ yields all conformal 
diffeomorphisms of $\S^n$.

For a conformal immersion $f:M\to \S^n$ of Riemann surface $M$, 
we may associate a \emph{(local) lift} to it, which is 
just a map $F$ from $M$ into the light cone $\L$ such that
the null line spanned by $F(p)$ is $f(p)$, $p\in M$. 
Taking derivatives with respect to a given complex coordinate $z$,
we find $\<F_z,F_{\zb}\> >0$ since $f$ is immersed, and
$\<F_z,F_z\>=0$ since $f$ is conformal. Furthermore, there is a 
rank-$4$ subbundle of $M\times\R^{n+1,1}$:
\[
V={\rm span} \{ F,\Re (F_z),\Im (F_z), F_{z\zb}, \}
\]
which is independent to the choice of local lift $F$ 
and complex coordinate $z$, hence well-defined. 
With signature $(3,1)$, $V$ describes
a M{\"o}bius invariant 2-sphere congruence, called the 
\emph{mean curvature sphere} of $f$. This name comes from 
the property that it is tangent to the surface and has the 
same mean curvature as the surface at the tangent point 
when the ambient space is endowed with a metric of Euclidean 
space (or any space form). 

Fix a local coordinate $z$. We find the canonical 
lift $Y$ determined by $\abs{\d Y}^2=\abs{\d z}^2.$ 
A M{\"o}bius invariant frame of $V\otimes\C$ is given as 
\begin{equation}
\label{frame}
\{Y,Y_z,Y_{\zb},N\},
\end{equation}
where we choose $N\in \Gamma(V)$ with
$\< N,N\>=0,~\< N,Y\>=-1,~\< N,Y_z\> =0.$ 
These frame vectors are orthogonal to each other except that
$\< Y_{z},Y_{\zb} \> = \frac{1}{2},~\< Y,N \> = -1.$

The fundamental equation in our surface theory is
\[
Y_{zz} + \frac{s}{2} Y = \kappa,
\]
which holds true for some complex valued function $s$
and a section $\kappa\in \Gamma(V^{\bot}\otimes\C)$,
because we may decompose $Y_{zz}$ according to 
$M\times\R^{n+1,1}=V\oplus V^{\bot}$, and $Y_{zz}$ is 
orthogonal to $Y,Y_z,Y_{\zb}$. In this way we introduce 
two basic M{\"o}bius invariants associated with 
the surface $f:M\to \S^n$: $s$ is interpreted as the 
Schwarzian derivative of immersion $f$, 
and $\kappa$ may be identified with the normal-valued 
Hopf differential up to scaling
(so $\kappa$ vanishes exactly at the umbilic points). 
The \emph{Schwarzian} $s$ and the \emph{Hopfian} $\kappa$ 
form a complete system of conformal invariants for surfaces 
in $\S^3$. Later we will need the fact that
$\kappa\frac{\d z^2}{\abs{\d z}}$
is a globally defined form. For more details see \cite{bpp}.

Let $\xi\in\Gamma(V^{\bot})$ be an arbitrary section of 
the normal bundle $V^{\bot}$, $D$ be the normal connection. 
We are now able to write down the structure equations:
\begin{equation} 
\label{mov-eq}
\left\{
\begin{array}{llll}
Y_{zz} &=& -\frac{s}{2} Y + \kappa, \\[.2cm]
Y_{z\zb} &=& 
-\<\kappa,\bar{\kappa}\> Y + \frac{1}{2}N, \\[.2cm]
N_{z} &=& -2 \<\kappa,\bar{\kappa}\> Y_z 
- s Y_{\zb}+ 2D_{\zb}\kappa, \\[.2cm]
\xi_z &=& D_z\xi + 2\<\xi,D_{\zb}\kappa\> Y
- 2\<\xi,\kappa\> Y_{\zb}.
\end{array}
\right.
\end{equation}
The conformal Gauss, Codazzi and Ricci equations 
as integrable conditions are 
\begin{eqnarray}
& \frac{1}{2}s_{\zb} = 3\< D_z\bar{\kappa},\kappa\>
+ \<\bar{\kappa},D_z\kappa\>, \label{gauss}\\[.1cm]
& \Im ( D_{\zb}D_{\zb}\kappa
+ \frac{\bar{s}}{2}\kappa ) = 0, \label{codazzi}\\[.1cm]
& R_{\zb z}^D\xi
:= D_{\zb}D_z\xi - D_z D_{\zb}\xi
= 2\<\xi,\kappa\>\bar{\kappa}
- 2\<\xi,\bar{\kappa}\>\kappa. \label{ricci}
\end{eqnarray}

Since $\kappa\frac{\d z^2}{\abs{\d z}}$ is independent to
the choice of coordinate $z$, there is a well-defined metric
\[
g:=e^{2\omega}\abs{\d z}^2
=4\<\kappa,\bar\kappa\>\abs{\d z}^2
\] 
over $M$, which is invariant under M\"obius transformations
and called the \emph{M\"obius metric}. 
Its Laplacian and Gaussian curvature $\K$ (called the 
\emph{M\"obius curvature}) are
\[
\Delta:=4e^{-2\omega}\partial_z\partial_{\bar z},\quad
\K:=-\Delta\omega=-4e^{-2\omega}\omega_{z\zb}.
\]
It is well known that
this metric is induced from the mean curvature sphere 
congruence and is conformal to the induced metric of $f$.
Hence the mean curvature sphere is also known as the
\emph{conformal Gauss map}\cite{bryant2}. The area of 
$M$ with respect to the M\"obius metric  
\[
W(f):=2{\rm i}\cdot\int_M \abs{\kappa}^2 \d z\wedge\d\zb
\]
is called \emph{Willmore functional} in this paper.
It coincides with the usual definition
\[
\widetilde{W}(f):=\int_M (H^2-K)\d M
\]
when $f$ is an immersed surface
in $R^3$ with mean curvature $H$ and Gauss curvature $K$.
The critical surfaces with respect to $W(f)$ are called
\emph{Willmore surface}, which we are most
interested in. They are characterized by the 
\emph{Willmore equation} \cite{bpp}:
\begin{equation}
\label{willmore}
D_{\zb}D_{\zb}\kappa+\frac{\bar{s}}{2}\kappa = 0.
\end{equation}
Note that this is stronger than the Codazzi equation \eqref{codazzi}.

\begin{remark}
\label{rem-associated}
When the data $\{s,\kappa\}$ satisfy
\eqref{gauss}\eqref{willmore}\eqref{ricci}, the data 
$\{s,e^{\rm it}\kappa\}$ also satisfy these integrable equations, 
where $t$ is a real parameter. This yields the \emph{associated 
family of Willmore surfaces}. 
\end{remark}

In the special case that $D_{\zb}\kappa$ depends linearly on
$\kappa$, locally there is a function $\mu$ such that
\begin{equation}
\label{swillmore}
D_{\zb}\kappa +\frac{\bar\mu}{2}\kappa =0.
\end{equation}
This time the surface is Willmore iff
\begin{equation}
\label{willmore2}
\mu_z - \frac{1}{2}{\mu}^2 - s= 0.
\end{equation}
Such a Willmore surface is called \emph{S-Willmore}
surface, which might be characterized by the property that 
it allows a dual Willmore surface envelopping the same mean 
curvature spheres \cite{ma}. Examples include Willmore 
surfaces in $\S^3$, isotropic surfaces in $\S^4$, 
and minimal surfaces in $\R^n, \S^n, \H^n$. The first surface class
is S-Willmore since the codimension equals one. The second class
has isotropic $\kappa$, whose derivative has to be
linearly dependent on itself due to codim-2 (we can verify
easily that it is also Willmore). Any minimal surface in a 
real space form has a dual surface as well. In general, 
the dual surface might be written down explicitly as
\begin{equation}
\label{yhat}
\widehat{Y}=\frac{1}{2}\abs{\mu}^2 Y 
+\bar\mu Y_z +\mu Y_{\zb} +N
\end{equation}
with respect to the frame $\{Y,Y_z,Y_{\zb},N\}$.
Simple calculation using \eqref{mov-eq} and 
\eqref{swillmore}\eqref{willmore2} yields
\begin{equation}
\label{yhat-z}
\widehat{Y}_z=\frac{\mu}{2}\widehat{Y}
   + \rho\left(Y_z + \frac{\mu}{2}Y\right),~~
   \text{where}~\rho:=\bar\mu_z-2\<\kappa,\bar\kappa\>.
\end{equation}
Here $\rho\abs{\d z}^2$ is a globally defined invariant 
associated with the surface pair $Y,\widehat{Y}$ \cite{ma}.
There follows
\begin{equation}
\label{yhat-metric}
\<\widehat{Y}_z,\widehat{Y}_z\>=0, \quad
\<\widehat{Y}_z,\widehat{Y}_{\zb}\>
=\frac{1}{2}\abs{\rho}^2.
\end{equation}
It is straightforward to verify that $\widehat{Y}$ share
the same mean curvature sphere as $Y$. So $\widehat{Y}$
determines a conformal Willmore immersion into $\S^n$ 
when $\rho\ne 0$ (see \cite{ma} and references therein).
\begin{lemma}
\label{lem-rho}
For S-Willmore surfaces with $\mu$ and $\rho$ defined as 
before, we have $\rho_{\zb}=\bar\mu\rho$. There exist a 
4-form $\Theta$ and a function $\P$ globally defined over $M$:
\[
\Theta\cdot \d z^4:=\rho\<\kappa,\kappa\>\d z^4,
\quad  \P:=\frac{\rho}{\<\kappa,\bar\kappa\>}.
\]
Moreover, $\Theta$ is holomorphic.
\end{lemma}
\begin{proof}
The first equality comes from \eqref{willmore2} and
\eqref{gauss}. The holomorphicity of $\Theta$
follows easily. Since 
$\rho\abs{\d z}^2$ and $\kappa\frac{\d z^2}{\abs{\d z}}$
are globally defined, the same is true for $\Theta$ and $\P$.
\end{proof}

The holomorphic 4-form $\Theta$ was first introduced 
in \cite{bryant2} for Willmore surfaces in $\S^3$. 
The function $\P$ will play a central role in the 
later discussions.

\section{Willmore surfaces with $\K\equiv c$ in $\S^3$}

In this section, $y:M\to \S^3$ is always assumed to be a 
Willmore surface without umbilic points. Codim-1 enables us
to choose a unit normal vector $X$ canonically,
which is exactly the conformal Gauss map from $M$ into
the de Sitter sphere $\S^{3,1}\subset \R^{4,1}$. The mean 
curvature sphere $\mathbb{P}(\L\cap V)$ is determined by
$V={\rm span} \{ Y,\Re (Y_z),\Im (Y_z), N, \}$
or by its orthogonal complement $U=V^{\bot}={\rm span}\{X\}$.
In this way we identified the conformal Gauss map and the
mean curvature sphere of $y$.

Now the normal-valued Hopf differential can be written as 
$\kappa\cdot X$, where we regard $\kappa$ as 
a scalar function. Denote the complex valued function $\P$ as
\[
\P=\abs{\P}\cdot e^{\rm i\psi}
=\frac{\rho}{\kappa\bar\kappa}.
\]
Between scalar invariants $\K$ and $\P$ exists important
connections as below:
\begin{lemma}
\label{lem-P}
\quad\\[-1.3cm]
\begin{subequations} 
\begin{align} 
&\Re(\P) =2(\K-1), \label{P1}\\
\text{and}~~~\Delta &(\log \P )=4\K+{\rm i}\cdot\Im (\P ),
\quad\text{when}\abs{\P}\ne 0. \label{P2}\\
\text{\hspace{-2cm} In particular,}
\quad\Delta &(\log\abs{\P})=4\K, \label{P3}\\
\Delta &\psi =\Im (\P )=\abs{\P}\sin\psi. \label{P4}
\end{align} 
\end{subequations}
\end{lemma}
\begin{proof}
From the definition of the M\"obius metric and function $\P$,  
\[
\K=-4e^{-2\omega}\omega_{z\zb}
=-\frac{1}{\kappa\bar\kappa}
\left [\frac{1}{2}\log (4\kappa\bar\kappa)\right ]_{z\zb}
=\frac{1}{4\kappa\bar\kappa}(\bar\mu_z+\mu_{\zb})
=\frac{1}{2}\Re(\P)+1
\]
by \eqref{swillmore}\eqref{yhat-z}. This proves \eqref{P1},
which is indeed true for any surface in $\S^3$ without the
Willmore condition. Next, using $\rho_{\zb}=\bar\mu\rho$
(Lemma~\ref{lem-rho}), there will be
\begin{align*}
\Delta(\log\P) &=4e^{-2\omega}
[\log\rho-\log (\<\kappa,\bar\kappa\> )]_{z\bar z}
=4e^{-2\omega}(\bar\mu_z-2\omega_{z\zb})\\
&=\frac{1}{\kappa\bar\kappa}(\rho+2\kappa\bar\kappa )+2\K 
=\P+2+2\K. 
\end{align*}
Note that this is valid for any S-Willmore surface in $\S^n$. 
Decomposing the right hand side into real and imaginary 
parts and invoking \eqref{P1} yields \eqref{P2}.
\end{proof}

The system \eqref{P1}\eqref{P2} has been obtained essentially 
in \cite{palmer} and \cite{li2}. After posing the restriction 
that the metric is of constant $\K$, we obtain an over-determined 
system, which might be solved completely in the sequel. 
The first result is a well-known fact.
\begin{proposition}
\label{prop-minimal}
Let $y:M\to \S^3$ be a Willmore surface. 
Then the following three conditions are equivalent:
1) $\P\equiv 0$;
2) $\K\equiv 1$;
3) $y$ is the image of a minimal surface in $\R^3$ 
under the inverse of a stereographic projection.
\end{proposition}
\begin{proof}
We assert 3)$\Leftrightarrow$1). First $\P\equiv 0$ 
implies $\rho\equiv 0$, which means $\widehat{Y}$ degenerates 
to a single point by \eqref{yhat-metric}. Regarding $y$
as in an affine space with $\widehat{Y}$ as the point 
at infinity, then every mean curvature sphere of $y$ 
turns out to be a plane, which is of mean curvature zero.
So $y$ itself is a minimal surface in this affine $\R^3$. 
The converse is true for the same reason.

Note 1)$\Rightarrow$2) by \eqref{P1}. So we need only to show 
2)$\Rightarrow$1). Suppose $\P\ne 0$ and $\K=1$ on an open subset. 
Then \eqref{P1} implies $\Re(F)=0$, hence $\psi=\arg\P=\pi/2$
is a constant. It follows from \eqref{P4} that $\P=0$, 
a contradiction.
\end{proof}

\begin{remark}
CMC-1 (constant mean curvature 1) surfaces in the 
hyperbolic space $H^3(-1)$ also satisfy $\K\equiv 1$. 
Such surfaces together with their minimal cousins 
are characterized as the only isothermic surfaces in $\S^3$ 
with $\K\equiv 1$ \cite{jeromin}.
\end{remark}

\begin{theorem}
\label{thm-clifford1}
Let $y:M\to \S^3$ be a Willmore surface. 
Then the following three conditions are equivalent:
1) $\abs{\P}=constant\ne 0$;
2) $\K\equiv 0$;
3) $y$ is locally M\"obius equivalent to the Clifford torus
with $\P=-2$.
\end{theorem}
\begin{proof}
As well-known Clifford torus is a homogeneous torus, thus 
$\K\equiv 0$, and $\P$ is a constant. Since it is Willmore, 
we may use \eqref{P1}\eqref{P2} to find $\P=-2$. 
So 3)$\Rightarrow$1),2).

Next we prove 1)$\Rightarrow$3).
$\abs{\P}=constant\ne 0$ implies $\K\equiv 0$ and $\Re(\P)=-2$ 
by \eqref{P1}\eqref{P2}. So $\psi=\arg \P=constant$.
It follows from \eqref{P4} that $\sin\psi=0$, hence $\psi=0$
and $\Im(\P)=0$. We conclude that $\P=-2$ and $\K=0$.
This time the M\"obius metric is flat, which enables us to find 
a suitable coordinate $z$ so that 
$\abs{\kappa}^2=\frac{1}{4}e^{2\omega}=1$. 
Consequently, $\rho=\P=-2$ is a constant, and $\mu=0$ due to 
$\rho_{\zb}=\bar\mu\rho$ (Lemma~\ref{lem-rho}).
For the basic invariants, 
$s=\mu_z-\frac{1}{2}\mu^2=0$ by \eqref{willmore2}, and 
$\kappa_{\zb}=D_{\zb}\kappa=0$ by \eqref{swillmore}.
Such a holomorphic and unitary function $\kappa$ must 
be a constant. After a proper rotation of the coordinate $z$
if necessary, we may assume $\kappa=1$ identically.
So all these invariants are constant now, and the frame
$\{{Y,Y_z,Y_{\zb}},N,X \}$ form a linear PDE system 
with constant coefficients:
\begin{displaymath}
\left\{
\begin{aligned} 
Y_{zz} &= ~X,\\
Y_{z\zb} &= -Y+\tfrac{1}{2}N,\\
N_z~ &= -2~Y_z, \\
X_z~ &= -2~Y_{\zb}.
\end{aligned} 
\right.
\end{displaymath}
It is obviously integrable with a unique solution
when the frame $\{{Y,Y_z,Y_{\zb}},N,X \}$ is given at
a fixed $p\in M$. Since the Clifford torus is a Willmore 
surface with $\K=0$, hence also a solution to this system, 
it must be congruent to $y$ up to a M\"obius transformation.
Indeed it may be written down explicitly at here:
\[
Y(u,v)=\frac{1}{2\sqrt{2}}\cdot \left (\sqrt{2},
\cos{2\sqrt{2}}u,\sin{2\sqrt{2}}u,
\cos{2\sqrt{2}}v,\sin{2\sqrt{2}}v \right ).
\]
Then one can verify easily that it solves the system above.

The most interesting part is to show that 2) implies 1) and 3). 
Suppose $\K\equiv 0$. By \eqref{P1}, 
$\Re(\P)=-2, \Im(\P)=\Re(\P)\cdot\tan\psi=-2\cdot\tan\psi.$
Thus $\P\ne 0$ and we may use \eqref{P4} to obtain
\begin{equation}
\label{tan}
\psi_{z\bar z}=-\frac{1}{2}e^{2\omega}\tan\psi.
\end{equation}
On the other hand, $\Delta(\log\abs{\P})=4\K=0$ 
by \eqref{P3}, so $\abs{\P}^2=4f\bar{f}$ 
for some non-zero holomorphic function $f$.
Invoking $\abs{\P}^2=(1+\tan^2\psi)\Re(\P)^2$, there is
\[
\tan\psi=\sqrt{f\bar{f}-1},\quad \text{i.e.}~~~
\psi=\arctan\sqrt{f\bar{f}-1}.
\]
Suppose $f\bar{f}>1$. Direct computation shows
\[
\psi_z = \frac{f_z}{f}\cdot\frac{1}{2\sqrt{f\bar{f}-1}},
\quad \text{and}~~~
\psi_{z\bar z} = \frac{f_z \bar{f}_{\bar z}}{-4(f\bar{f}-1)^{3/2}}.
\]
Substitute these expressions into \eqref{tan}, we find
\[
e^{2\omega}=\frac{1}{2}\frac{\abs{f_z}^2}{(1-f\bar{f})^2}.
\]
But the metric given above has constant gaussian curvature $-8$
when $f$ is holomorphic. A contradiction! The only possibility 
left is that $f\bar{f}\equiv 1$ and $\psi\equiv 0, \P\equiv -2$. 
This proves 1), hence also 3).
\end{proof}

Theorem~\ref{thm-clifford1} gives a nice characterization 
of Clifford torus. Such results are always interesting, 
especially when they might be related to the famous
Lawson's Conjecture and Willmore Conjecture. 
For more characterization theorems about Clifford torus, 
see \cite{kenmotsu2,sasaki2,vlachos1,vlachos2}. 
Below is another one using global assumptions.
(The reader who is interested in general surfaces with 
$\K=0$ might refer to \cite{sulanke}.)

\begin{theorem}
\label{thm-clifford2}
Let $y:M\to \S^3$ be a Willmore surface which is complete 
with respect to the M\"obius metric. If $\K\le 0$ on $M$, then 
$\K\equiv 0$ and $y$ is M\"obius equivalent to Clifford torus.
\end{theorem}

\begin{proof}
We introduce a new metric on $M$ as
\begin{equation}
\label{metric}
\tilde{g}=\sqrt{\abs{\P}}g=\sqrt{\abs{\P}}e^{2\omega}\abs{\d z}^2
=e^{2\tilde\omega}\abs{\d z}^2,
\end{equation}
where $\tilde\omega=\tfrac{1}{4}\log\abs{\P}+\omega$.
From the assumption $\K\le 0$ and \eqref{P4}, we have
$\Re(\P)\le -2$, hence $\abs{\P}\ge 2$. So the new metric
is positive definite and complete, too. 
The Gaussian curvature of $(M,\tilde{g})$ is
\[
\widetilde{K}=-4e^{-2\tilde\omega}(\tilde\omega)_{z\bar z}=0,
\]
by \eqref{P3}. We may assume that $M$ is simply connected.
Therefore $(M,\tilde{g})$ is an Euclidean plane,
and $(M,g)$ is conformally equivalent to the complex plane $\C$. 
Since $\abs{\P}\ge 2$ and $\K\le 0$, by \eqref{P3}, 
$\log\abs{\P}$ is a subharmonic function on $(\C,\tilde{g})$ 
with lower bound $\log 2$. We conclude that $\abs{\P}=constant>0$.
The proof is finished by using Theorem~\ref{thm-clifford1}.
\end{proof}

Now we can prove the first main theorem.

\begin{theorem}
\label{thm-constant1}
Let $y:M\to \S^3$ be a Willmore surface with constant M\"obius 
curvature $\K$. Then $y$ is M\"obius equivalent 
either to a minimal surface in $\R^3$ with $\K=1,\P=0$,
or to the Clifford torus with $\K=0,\P=-2$.
\end{theorem}
\begin{proof}
Let's assume $\K=constant\ne 0,1$. This time $\psi$ could not 
be a constant. (Otherwise it contradicts the assumption
by \eqref{P1}\eqref{P2}.) The first step is to show that 
this real function $\psi:M\to \R$ is an isoparametric function.

First, $\Re(\P)=2(\K -1)$ is a non-zero constant. Note
\begin{alignat}{2}
\Im(\P) &= \Re(\P)\cdot \tan\psi &&=2(\K -1)\tan\psi, 
\label{imp}\\
\abs{\P} &= \Re(\P)/\cos\psi &&= 2(\K -1)\sec\psi.
\label{modp}
\end{alignat}
Substitute \eqref{imp} and \eqref{modp} separately into 
$\Delta\psi=\Im(\P)$ \eqref{P4}
and $\Delta\log\abs{\P}=4\K$ \eqref{P3}, we find
\begin{align}
\Delta\psi~ &= 2(\K -1)\tan\psi=G(\psi), \label{isop1}\\
||\nabla\psi||^2 &= 4e^{-2\omega}\abs{\psi_z}^2 = 
4\left (\K\cos^2\psi-\frac{\K -1}{2}\sin^2\psi\right)=F(\psi).
\label{isop2}
\end{align}
Here \eqref{isop2} follows from \eqref{P3}\eqref{isop1}
via the following computation:
\begin{align*}
e^{2\omega}K &=(\log\abs{\P})_{z\bar z}=(-\log\cos\psi)_{z\bar z}
=\psi_{z\bar z}\cdot\tan\psi
+\psi_z \psi_{\bar z}\sec^2\psi \\
&=\frac{\K -1}{2}e^{2\omega}\tan^2\psi
+\psi_z \psi_{\bar z}\sec^2\psi.
\end{align*}
Such a function $\psi$ satisfying
$\Delta\psi=G(\psi),~||\nabla\psi||^2 =F(\psi)$ 
is called an \emph{isoparametric function}, 
where $F,G$ are two functions of $\psi$. 

We want to find contradiction in these equations.
This is done by appealing to the following identity
(see Lemma~\ref{lem-isop}):
\begin{equation}
\label{isop}
2\K F+(2G-F')(G-F')+F(2G'-F'')=0,
\end{equation}
where $G'=\d G/\d \psi,F'=\d F/\d \psi$, etc.. 
Substituting the expressions of $F(\psi)$ and $G(\psi)$
\eqref{isop1}\eqref{isop2} into this identity, 
after simplification we find
\[
4(27\K -8)(\K -1)-4(3\K -1)(3\K -8)\cos^2\psi =0.
\]
The left hand side is a polynomial of $\cos\psi$.
It is identically zero iff all the coefficients vanish.
Yet this is obviously impossible for any constant $\K$.
A contradiction! Thus our assumption in the beginning
is wrong. The conclusion follows from Proposition~
\ref{prop-minimal} and Theorem~\ref{thm-clifford1}.
\end{proof}

\begin{remark}
Gorokh\cite{gorokh2} studied minimal surfaces of constant 
Gaussian curvature in the 4-dim de Sitter's space $\S^{3,1}(1)$.
Note that such surfaces are equivalent to Willmore surfaces
of constant M\"obius curvature in $\S^3$. Any surface of the 
former class, denoted as $X$, is exactly the conformal Gauss map 
of one in the latter class, and vice versa.
Except the case $\K=1$ where $X$ is totally umbilic,
he obtained the same result as Theorem~\ref{thm-constant1}.
(He announced that the exceptional case $\K=1$, also the 
easiest case, would be indicated in his next paper, 
which did not appear.) His method is similar to our previous
proof in essence, yet buried in tedious computations. 
We think that the same result in terms of Willmore surfaces
might interest more people, and our presentation not only 
clarifies the key points, but also illustrates a 
general method for solving such problems.
\end{remark}

As a by-product of \eqref{P1}\eqref{P2}, we find the following
characterization of the associated Willmore surfaces.

\begin{proposition}
\label{prop-associated}
Let $y,\tilde{y}:M\to \S^3$ be two conformal Willmore immersions. 
The invariants of $y$ satisfy $\P\ne 0$ and $\K\ne 0$.
Then $\P=\widetilde{\P}$ iff 
$y$ and $\tilde{y}$ are in the same associated family.
\end{proposition}
\begin{proof}
When $y$ and $\tilde{y}$ are in the same associated family, 
their Hopfians are related by $\tilde\kappa=e^{{\rm i}t}\kappa$ 
for a constant $t\in \R$. Thus they share the same $\mu$
and $\<\kappa,\bar\kappa\>$ (M\"obius metric). By definition 
$\P=\widetilde{\P}$. So we need only to show the converse. 

Choose the same coordinate $z$ for $y$ ad $\tilde{y}$. 
Then it follows from \eqref{P1} that $\K=\widetilde{\K}$. 
By \eqref{P3}, $(e^{-2\omega}-e^{-2\tilde\omega})
(\log\abs{\P})_{z\bar z}=0$.
Because $\K\ne 0$, $(\log\abs{\P})_{z\bar z}\ne 0$, 
so $\omega=\tilde\omega$, i.e. $x,\tilde{x}$ induce 
the same M\"obius metric. Combined with 
$\P=\widetilde{\P}\ne 0$, we have $\rho=\tilde\rho\ne 0$. 
Taking logarithm at both sides and differentiating yields 
$\mu=\tilde\mu$ by Lemma~\ref{lem-rho}. As the consequence,
their Schwarzians are equal due to \eqref{willmore2},
and their Hopfians differ by a constant factor, i.e.
$\tilde\kappa=e^{{\rm i}t}\kappa$, due to \eqref{swillmore}.
So $y,\tilde{y}$ are in the same associated family
of Willmore surfaces.
\end{proof}

\section{Isotropic surfaces with $\K\equiv c$ in $\S^4$}

A surface in $\S^4$ is isotropic if $\<\kappa,\kappa\>=0$.
This is a M\"obius invariant notion, although usually
defined in terms of metric geometry. For its important 
geometric meaning we refer to \cite[Ch.~8]{burstall}.

The isotropic condition implies that all derivatives
of the Hopfian, $\{D_z\kappa,D_{\zb}\kappa,\cdots\}$,
are contained in the isotropic line bundle 
${\rm span}\{\kappa\}$. As a corollary of this fact
and the conformal Codazzi equation \eqref{codazzi},
such surfaces must be S-Willmore.

From now on we always assume that $y:M\to \S^4$ is 
an isotropic surface without umbilic points. 
So $\kappa\ne 0$ and we have
$D_{\zb}\kappa=-\frac{\bar\mu}{2}\kappa,
~D_z\kappa=\lambda\kappa$
for some locally defined function $\mu,\lambda$.
Substitute them into the conformal Ricci equation
\eqref{ricci}
\[
D_{\zb}D_z\kappa=D_z D_{\zb}\kappa 
-2\<\kappa,\bar\kappa\>\kappa
+2\<\kappa,\kappa\>\bar\kappa.
\]
It follows
$\bar\lambda_z=-\frac{\mu_{\zb}}{2}-2\<\kappa,\bar\kappa\>$.
Then the M\"obius curvature is computed as below:
\begin{align*}
\K &=-\frac{1}{2\<\kappa,\bar\kappa\>}
\big (\log\<\kappa,\bar\kappa\>\big )_{z\zb}
=-\frac{1}{2\<\kappa,\bar\kappa\>}
\left (-\frac{\bar\mu}{2}+\bar\lambda\right )_z \\
&=\frac{1}{4\<\kappa,\bar\kappa\>}
\big (\rho+\bar\rho+8\<\kappa,\bar\kappa\>\big )
=\frac 12\Re(\P)+2,
\end{align*}
where $\rho$ and $\P$ are defined as before.
On the other hand, in the proof to Lamma~\ref{lem-P},
it was noticed that for any S-Willmore surface in $\S^n$ holds
$\Delta(\log\P)=\P+2+2\K$
if $\P\ne 0$. Combined with the previous formula, 
we have proved 
\begin{lemma}
\label{lem-P2}
For any isotropic surface satisfying $\P\ne 0$, there are
\begin{subequations} 
\begin{align} 
&\K~ =\tfrac 12\Re (\P )+2=\tfrac 12\abs{\P}\cos\psi+2, 
\label{P5}\\
&\Delta (\log\abs{\P})=4\K -2, \label{P6}\\
&\Delta \psi =\Im (\P )=\abs{\P}\sin\psi, \label{P7}
\end{align} 
\end{subequations}
\end{lemma}
Where $\psi$ is the argument of $\P$.
Similar as in the previous section, with the help 
of this lemma we prove two characterization theorems, 
then complete the final classification.
\begin{theorem}
\label{thm-complex}
Let $y:M\to \S^4$ be an isotropic surface 
without umbilic points. 
Then the following three conditions are equivalent:
1) $\P\equiv 0$;
2) $\K\equiv 2$;
3) Locally $y$ is M\"obius equivalent to a complex curve in 
$\C^2\cong \R^4$, where the identification is made by 
choosing a suitable orthogonal complex structure on $\R^4$.
\end{theorem}
\begin{proof}
If $\P=0$ on an open subset of $M$. Then the dual Willmore 
surface $\widehat{Y}$ will degenerate to a single point as we 
observed before. This case $y$ is M\"obius equivalent to a
minimal surface in an affine 4-space. It is known that
an isotropic minimal surface in $\R^4$ is a complex curve
in $\R^4\cong\C^2$ . (See \cite{hoffman} for a proof. 
The straightforward way is to define a complex structure $J$ 
on $TM$ and $T^{\bot}M$, then show its parallel translation 
$\nabla J=0$. Hence it extends to a complex structure on $\R^4$.)
By \eqref{P5}, $\K\equiv 2$. Conversely, $\K\equiv 2$ implies
$\P\equiv 0$ like in Proposition~\ref{prop-minimal}. 
So complex curves in $\C^2$ are isotropic, minimal and $\K=2$,
and they are characterized by the 
combination of any two of these conditions.
\end{proof}

\begin{theorem}
\label{thm-veronese}
Let $y:M\to \S^4$ be an isotropic surface without 
umbilic points. 
Then the following three conditions are equivalent:
1) $\abs{\P}=constant \ne 0$;
2) $\K\equiv \frac 12$;
3) Locally $y$ is M\"obius equivalent to the Veronese sphere 
in $\S^4$ with $\P\equiv -3$.
\end{theorem}
\begin{proof}
The Veronese 2-sphere is also a homogeneous 
Willmore surface, thus $\K, \P$ are both constant, 
which could be found as $\K=\frac 12,\P=-2$ 
by \eqref{P5}-\eqref{P7}. (The case $\P=0$ is impossible
because we know that the Veronese surface is only minimal
in $\S^4$ and not minimal in any Euclidean space. Indeed
$\K=\frac 12$ was a known fact about this surface.) 
So 3)$\Rightarrow$1),2).

Next we prove 1)$\Rightarrow$2),3). $\abs{\P}=constant\ne 0$ 
implies $\K=\frac 12$ and $\Re(\P)=-3$ by \eqref{P6}\eqref{P5}. 
So the argument $\psi$ is constant.
It follows from \eqref{P7} that $\sin\psi=0$, hence $\psi=0$
and $\Im(\P)=0, \P=-3$. By the definitions of $\P$
and Lemma~\ref{lem-rho}, one finds
\[
\rho=-3\<\kappa,\bar\kappa\>,\quad
\bar\mu=\big (\log\rho\big )_{\zb}
=\big (\log\<\kappa,\bar\kappa\>\big )_{\zb}=2\omega_{\zb}.
\]
Together with $(\log\<\kappa,\bar\kappa\>)_{\zb}
=-\frac{\bar\mu}{2}+\bar\lambda$ there is 
$\bar\lambda=3\omega_{\zb}$, and the Schwarzian is
\[
s=\mu_z-\frac 12\mu^2=2(\omega_{zz}-\omega_z^2),
\]
which is holomorphic by \eqref{gauss} or by 
$\omega_{z\zb}=-\frac 18 e^{2\omega}$. 
Sum together, we have shown that the complex frame
$\{Y,Y_z,N,\kappa\}$ of this isotropic surface satisfies
\begin{displaymath}
\left\{
\begin{aligned} 
Y_{zz} &= (\omega _z^2 -\omega _{zz})Y+\kappa,\\
Y_{z\zb} &= -\tfrac 14 e^{2\omega}~Y+\tfrac 12 N,\\
N_z~ &= -\tfrac 12 e^{2\omega}~Y_z 
+2(\omega _z^2 -\omega _{zz})Y_{\zb}-2\omega _{\zb}\kappa, \\
\kappa _z~ &= D_z\kappa =3\omega _z\kappa,\\
\kappa _{\zb}~ &= D_{\zb}\kappa +2\<\kappa ,D_z\bar\kappa\>Y
-2\<\kappa,\bar\kappa\> Y_z \\
&= -\omega_{\zb}\kappa -\tfrac 12\omega_z e^{2\omega}Y
-\tfrac 12 e^{2\omega}Y_z
\end{aligned} 
\right.
\end{displaymath}
with respect to a given coordinate $z$ and 
a metric $e^{2\omega}\abs{\d z}^2$ of constant
curvature $\K=\frac 12$. This PDE system is integrable
by our construction and verification. The solution
is unique when the initial values of these frame vectors
at a given $p\in M$ are fixed. So $Y$ is locally congruent to 
the known solution, the Veronese sphere, 
up to a M\"obius transformation.

The final part ``2)$\Rightarrow$1),3)" will be included
in the proof to the next theorem.
\end{proof}

\begin{theorem}
\label{thm-constant2}
Let $y:M\to \S^4$ be an isotropic surface without umbilic 
points and $\K=constant$. Then locally $y$ is M\"obius 
equivalent either to a complex curve in $\C^2$ with $\K=2,\P=0$, 
or to the Veronese sphere with $\K=\frac 12,\P=-3$.
\end{theorem}
\begin{proof}
Suppose the constant $\K\ne 2$, or equivalently $\P\ne 0$,
on an open subset of $M$. Then $\Re(\P)=2\K -4$ is a 
non-zero constant. Put 
$\Im(\P)=2(\K -2)\tan\psi$ and $\abs{\P}=2(\K -2)\sec\psi$ into 
\eqref{P7}\eqref{P6}, we find
\begin{align}
\Delta\psi~ &= 2(\K -2)\tan\psi =G(\psi ), \label{isop3}\\
||\nabla\psi||^2 &=
4\left (\K\cos^2\psi-\frac{\K -2}{2}\sin^2\psi\right )=F(\psi )
\label{isop4}
\end{align}
just as we did in the proof to Theorem~\ref{thm-constant1}.
Substitute the expressions of $F(\psi)$ and $G(\psi)$
\eqref{isop3}\eqref{isop4} into the same identity 
\eqref{isop} for isoparametric function $\psi$:
\[
2\K F+(2G-F')(G-F')+F(2G'-F'')=0.
\]
We conclude in the same manner that when $\psi$ is a 
non-constant function there exists no such a $\K$ 
satisfying this identity. The left possibility is
$\psi=constant$. (In particular, the condition 
``$\K\equiv \frac 12$" in Theorem~\ref{thm-veronese} implies
$\psi=constant$.) This constant $\psi$ has to be zero
due to \eqref{P7}, hence $\P=\Re(\P)$ is also constant
by \eqref{P5}, and $\K\equiv \frac 12$. This surface 
is M\"obius equivalent to the Veronese 2-sphere by the
``1)$\Rightarrow$3)" part in Theorem~\ref{thm-veronese}.
Thus we complete the proof to Theorem~\ref{thm-veronese}
and Theorem~\ref{thm-constant2} at the same time.
\end{proof}
\begin{remark}
Like that for Clifford torus, characterization of 
Veronese surface(s) is also a favorite topic
in the study of submanifolds with abundant results.
Here we would like to emphasize that our characterizations
of Clifford torus and Veronese 2-sphere are in terms
of conformal invariants. Another such result is
\cite{li1}. 
\end{remark}

\section{Further remarks}

We note that our results are good examples of
the phenomenon of \emph{quantization} and 
\emph{non-negativity} in such problems. Namely, certain
geometric quantity (like the scalar curvature) of
those rigid objects (like minimal submanifolds), 
when being constant, can take
its value only in a discrete subset of $\R$. 
In particular, the negative values are usually 
forbidden, like minimal surfaces in $\S^n$ of
constant Gaussian curvature mentioned at the beginning.

The equation systems we met in Lemma~\ref{lem-P} and
Lemma~\ref{lem-P2} are examples of Toda system, an
important class of integrable systems. They are 
usually associated with certain surface
classes: besides Willmore surfaces in $\S^3$ and 
isotropic surfaces in $\S^4$, there are Willmore
surfaces in $\S^{2,1}$ \cite{alias}, minimal surfaces
in $\S^4$ \cite{kenmotsu2} and (non-super) minimal surfaces
in $\C\mathbb{P}^2$ \cite{chi} among others.

In such systems, under the assumption of constant 
Gaussian curvature of the metric in consideration,
one usually find some invariants being isoparametric 
functions. For results concerning such functions
the reader might consult \cite{chi} and references
therein. The most useful property of such functions
on a Riemannian surface is the identity \eqref{isop}.
Many people used it in their work \cite{chi,gorokh2,
kenmotsu2,kenmotsu3,sasaki2,vlachos1,vlachos2}. 
It is usually refered to Eisenhart's book 
\cite[p.164]{eisenhart}. Yet this identity is not given
explicitly there. Here we state it as the lemma below
with an independent proof.
\begin{lemma}
\label{lem-isop}
Let $(M,\d s^2)$ be a 2-dim Riemannian manifold and
$\psi:M\to\R$ be a smooth isoparametric function, i.e.
\[
||\nabla\psi||^2=F(\psi),\quad \Delta\psi=G(\psi)
\]
for smooth functions $F,G:\R\to\R$, where $\nabla$
and $\Delta$ denote the gradient and Laplacian with
respect to the metric $\d s^2$. 
Then on $\{p\in M: \nabla\psi(p)\ne 0\}$, 
the Gaussian curvature $\K$ satisfies the identity
\[
2\K F+(2G-F')(G-F')+F(2G'-F'')=0.
\]
Here $G'=\d G/\d \psi,F'=\d F/\d \psi$, etc..
\end{lemma}
\begin{proof}
Write $\d s^2=e^{2\omega}\abs{\d z}^2$. Then
$F(\psi)=4e^{-2\omega}\psi_z \psi_{\zb},
G(\psi)=4e^{-2\omega}\psi_{z\zb}$. Eliminating
$e^{2\omega}$, we find
\begin{equation}
\label{a1}
\psi_{z\zb}=f_z f_{\zb}\cdot\frac GF.
\end{equation}
On the other hand, there follows 
$e^{2\omega}=\abs{\frac{2\psi_z}{\sqrt{F}}}^2$. To compute
its Gaussian curvature, we notice that
\begin{align*}
\left (\frac{\psi_z}{\sqrt{F}}\right )_{\zb}
&= \frac{\psi_{z\zb}}{\sqrt{F}}
-\frac 12\frac{\psi_z}{F\sqrt{F}}F_{\zb}
=\frac{\psi_z \psi_{\zb}}{F\sqrt{F}}G
-\frac 12\frac{\psi_z}{F\sqrt{F}}\psi_{\zb}\cdot F'\\
&=\frac{\psi_z}{\sqrt{F}}\psi_{\zb}\cdot \frac{\d f}{\d \psi}
=\frac{\psi_z}{\sqrt{F}}\cdot f_{\zb}
\end{align*}
by \eqref{a1}, where $f(\psi)$ is a real function 
solving the ODE $f'=\frac 1F \left(G-\frac 12 F'\right)$.
There follows 
$\big(\log\dfrac{2\psi_z}{\sqrt{F}}-f\big)_{\zb}=0$,
hence $\dfrac{2\psi_z}{\sqrt{F}}=e^f\cdot h$ for
some holomorphic function $h$ on $M$. 
Write $e^{2\omega}=e^{2f}\cdot\abs{h}^2$. 
The curvature $\K$ satisfies
\begin{align*}
-\K &=4e^{-2\omega}\omega_{z\zb}=4e^{-2\omega}f_{z\zb}
=4e^{-2\omega}\left[\frac{\psi _{\zb}}{F} 
\big (G-\frac 12 F'\big )\right ]_z \\
&=4e^{-2\omega}\left[ \big (\frac{\psi_{z\zb}}{F} 
-\frac{\psi _z\psi_{\zb}}{F^2}\cdot F' \big )
\big (G-\frac 12 F'\big )+\frac{\psi_{\zb}}{F} \cdot\psi _z
\cdot\big (G'-\frac 12 F''\big )\right ] \\
&=\frac 1F \left[\big (G-F'\big )\big (G-\frac 12 F'\big )
+F\big (G'-\frac 12 F''\big )\right ].
\qedhere
\end{align*}
\end{proof}

\vspace{0.4cm}

\begin{tabbing}
xxxxxxxxxxxxxxxxxxxxxxxxxxxxxxxxxxxx\=\kill
Xiang Ma \> Changping Wang\\
LAMA \> LAMA\\
School of Mathematical Sciences \> School of Mathematical Sciences\\
Peking University \> Peking University\\
Beijing 100871 \> Beijing 100871\\
People's Republic of China \> People's Republic of China\\
Email: \textsf{maxiang@math.pku.edu.cn} \> Email: \textsf{wangcp@pku.edu.cn}
\end{tabbing}


\begin{thebibliography}{99}

\bibitem{alias}
Alias, L., Palmer, B.:
Conformal geometry of surfaces in Lorentzian space forms,
Geom. Dedicata {\bf 60}, 301--315(1996)

\bibitem{blaschke}
Blaschke, W.: Vorlesungen {\"u}ber Differentialgeometrie III: 
Differentialgeometrie der Kreise und Kugeln, 
Springer Grundlehren XXIX, Berlin, 1929

\bibitem{bryant1}
Bryant, R.: A duality theorem for Willmore surfaces, 
J. Diff. Geom. {\bf 20}, 23--53(1984)

\bibitem{bryant2}
Bryant, R.: 
Minimal surfaces of constant curvature in $S^n$, 
Trans. Amer. Math. Soc. {\bf 290}(1), 259--271(1985)

\bibitem{burstall}
Burstall, F., Ferus, D., Leschke, K., Pedit, F.,
Pinkall, U.: 
Conformal geometry of surfaces in $S^4$ and quaternions, 
Lecture Notes in Mathematics 1772, Springer, Berlin, 2002 

\bibitem{bpp}
Burstall, F., Pedit, F., Pinkall, U.: 
Schwarzian derivatives and flows of surfaces, 
Contemporary Mathematics {\bf 308}, 2002, pp.39--61

\bibitem{calabi}
Calabi, E.: 
Minimal immersions of surfaces in Euclidean spheres, 
J. Diff. Geom. {\bf 1}, 111--125(1967)

\bibitem{chi}
Chi, Q., Jensen, G., Liao, R.:
Isoparametric functions and flat minimal tori in $\C P^2$,
Proc. Amer Math. Soc. {\bf 123}(9), 2849--2854(1995)

\bibitem{eisenhart}
Eisenhart, L.:
An introduction to differential geometry, Princeton, NJ, 1947

\bibitem{fujiki}
Fujiki, M.:
On the Gauss map of minimal surfaces immersed in $R\sp n$.
Kodai Math. J. {\bf 9}(1), 44--49(1986)

\bibitem{gorokh1}
Gorokh, V.:
Two-dimensional minimal surfaces in a 
pseudo-Euclidean space,
Ukrain. Geom. Sb. {\bf 31}, 36--47(1988); translation in 
J. Soviet Math. {\bf 54}(1), 691--699(1991)

\bibitem{gorokh2}
Gorokh, V.: 
Minimal surfaces of constant Gaussian curvature in a
pseudo-Riemannian sphere, 
Ukrain. Geom. Sb. {\bf 32}, 27--34(1989); translation in
J. Soviet Math. {\bf 59}(2), 709--714(1992)

\bibitem{hoffman}
Hoffman, D., Osserman, R.:
The geometry of the generalized Gauss map,
Mem. Amer. Math. Soc. {\bf 236} 1980

\bibitem{jeromin}
Hertrich-Jeromin, U., Musso, E., Nicolodi, L.: 
M{\"o}bius geometry of surfaces
of constant mean curvature $1$ in hyperbolic space, 
Ann. Global Anal. Geom. {\bf 19}(2), 185--205(2001)

\bibitem{kenmotsu1}
Kenmotsu, K.: 
On minimal immersions of $R^2$ into $S^n$, 
J. Math. Soc. Japan {\bf 28}, 182-191(1976) 

\bibitem{kenmotsu2}
Kenmotsu, K.: 
Minimal surfaces with constant curvature in 
4-dimensional space forms, 
Proc. Amer. Math. Soc. {\bf 89}(1), 133--138(1983)

\bibitem{kenmotsu3}
Kenmotsu, K., Masuda, K.: 
On minimal surfaces of constant curvature in 
two-dimensional complex space forms,
J. reine angew. Math. {\bf 523}, 69--101(2000)

\bibitem{lawson}
Lawson, H. Blaine, Jr.:
The Riemannian geometry of holomorphic curves,
Bol. Soc. Brasil. Mat. {\bf 2}(1), 45--62(1971)


\bibitem{li1}
Li, H., Wang, C., Wu, F.:
A Moebius characterization of Veronese surfaces in $S^n$,
Math. Ann. {\bf 319}, 707--714(2001)

\bibitem{li2}
Li, H., Wang, C., Zhao, G.: 
A new M\"obius invariant function for surfaces in $S^3$, 
preprint, 2003

\bibitem{ma}
Ma, X.: 
Isothermic and S-Willmore Surfaces as Solutions to 
a Problem of Blaschke,
Results Math. {\bf 48}(3-4), 301--309(2005)

\bibitem{masaltsev}
Masaltsev, L.:
Minimal surfaces in $R\sp{5}$ whose Gaussian image 
has a constant curvature,
Mat. Zametki {\bf 35}(6), 927--932(1984)

\bibitem{nie}
Nie, C., Ma, X., Wang, C.:
Conformal CMC-surfaces in Lorentzian space forms,
Chin. Ann. Math. series B, to appear.

\bibitem{palmer}
Palmer, B.: 
The conformal Gauss map and the stability of Willmore surfaces,
Ann. Global Anal. Geom. {\bf 9}(3), 305--317(1991)

\bibitem{sasaki1}
Sasaki, M.:
Exceptional minimal surfaces whose Gauss images 
have constant curvature,
Tokyo J. Math. {\bf 15}, 381--388(1992)

\bibitem{sasaki2}
Sasaki, M.:
Minimal 2-tori in $S^4$ whose Gauss images 
have constant curvature,
Arch. Math. {\bf 62}(5), 470--474(1994)

\bibitem{sasaki3}
Sakaki, M.:
A note on minimal surfaces in $\mathbf R\sp n$ 
whose Gauss images have constant curvature.
Bull. Fac. Sci. Technol. Hirosaki Univ. 
{\bf 3}(2), 79--82(2001)

\bibitem{sulanke}
Sulanke, R.:
M\"obius geometry. VII. On channel surfaces. 
Proceedings of the 3rd Congress of Geometry 
(Thessaloniki, 1991), 410--419, 
Aristotle Univ. Thessaloniki, Thessaloniki, 1992.

\bibitem{vlachos1}
Vlachos, T.:
The third fundamental form of minimal surfaces in a sphere,
Arch. Math. {\bf 74}, 66--74(2000)

\bibitem{vlachos2}
Vlachos, T.:
A characterization of the Clifford torus,
Arch. Math. {\bf 85}, 175--182(2005)

\bibitem{wang}
Wang, C.: 
Moebius geometry of submanifolds in $S^n$, 
Manuscripta Math, {\bf 96}, 517--534(1998)


\end{thebibliography}
\end{document}